\documentclass[12pt]{amsart}
\usepackage{amssymb}
\usepackage{amsmath, amsfonts}

\usepackage[mathscr]{eucal}

\input xypic
\xyoption{all}

\makeindex
\makeglossary

\begin{document}

\baselineskip = 16pt

\newcommand \ZZ {{\mathbb Z}}
\newcommand \NN {{\mathbb N}}
\newcommand \RR {{\mathbb R}}
\newcommand \PR {{\mathbb P}}
\newcommand \AF {{\mathbb A}}
\newcommand \GG {{\mathbb G}}
\newcommand \QQ {{\mathbb Q}}
\newcommand \CC {{\mathbb C}}
\newcommand \bcA {{\mathscr A}}
\newcommand \bcC {{\mathscr C}}
\newcommand \bcD {{\mathscr D}}
\newcommand \bcF {{\mathscr F}}
\newcommand \bcG {{\mathscr G}}
\newcommand \bcH {{\mathscr H}}
\newcommand \bcM {{\mathscr M}}
\newcommand \bcI {{\mathscr I}}
\newcommand \bcJ {{\mathscr J}}
\newcommand \bcK {{\mathscr K}}
\newcommand \bcL {{\mathscr L}}
\newcommand \bcO {{\mathscr O}}
\newcommand \bcP {{\mathscr P}}
\newcommand \bcQ {{\mathscr Q}}
\newcommand \bcR {{\mathscr R}}
\newcommand \bcS {{\mathscr S}}
\newcommand \bcT {{\mathscr T}}
\newcommand \bcV {{\mathscr V}}
\newcommand \bcU {{\mathscr U}}
\newcommand \bcW {{\mathscr W}}
\newcommand \bcX {{\mathscr X}}
\newcommand \bcY {{\mathscr Y}}
\newcommand \bcZ {{\mathscr Z}}
\newcommand \goa {{\mathfrak a}}
\newcommand \gob {{\mathfrak b}}
\newcommand \goc {{\mathfrak c}}
\newcommand \gom {{\mathfrak m}}
\newcommand \gon {{\mathfrak n}}
\newcommand \gop {{\mathfrak p}}
\newcommand \goq {{\mathfrak q}}
\newcommand \goQ {{\mathfrak Q}}
\newcommand \goP {{\mathfrak P}}
\newcommand \goM {{\mathfrak M}}
\newcommand \goN {{\mathfrak N}}
\newcommand \uno {{\mathbbm 1}}
\newcommand \Le {{\mathbbm L}}
\newcommand \Spec {{\rm {Spec}}}
\newcommand \Gr {{\rm {Gr}}}
\newcommand \Pic {{\rm {Pic}}}
\newcommand \Jac {{{J}}}
\newcommand \Alb {{\rm {Alb}}}
\newcommand \alb {{\rm {alb}}}
\newcommand \Corr {{Corr}}
\newcommand \Chow {{\mathscr C}}
\newcommand \Sym {{\rm {Sym}}}
\newcommand \Prym {{\rm {Prym}}}
\newcommand \cha {{\rm {char}}}
\newcommand \eff {{\rm {eff}}}
\newcommand \tr {{\rm {tr}}}
\newcommand \Tr {{\rm {Tr}}}
\newcommand \pr {{\rm {pr}}}
\newcommand \ev {{\it {ev}}}
\newcommand \cl {{\rm {cl}}}
\newcommand \interior {{\rm {Int}}}
\newcommand \sep {{\rm {sep}}}
\newcommand \td {{\rm {tdeg}}}
\newcommand \alg {{\rm {alg}}}
\newcommand \im {{\rm im}}
\newcommand \gr {{\rm {gr}}}
\newcommand \op {{\rm op}}
\newcommand \Hom {{\rm Hom}}
\newcommand \Hilb {{\rm Hilb}}
\newcommand \Supp{{\rm Supp}}
\newcommand \Sch {{\mathscr S\! }{\it ch}}
\newcommand \cHilb {{\mathscr H\! }{\it ilb}}
\newcommand \cHom {{\mathscr H\! }{\it om}}
\newcommand \colim {{{\rm colim}\, }}
\newcommand \End {{\rm {End}}}
\newcommand \coker {{\rm {coker}}}
\newcommand \id {{\rm {id}}}
\newcommand \van {{\rm {van}}}
\newcommand \spc {{\rm {sp}}}
\newcommand \Ob {{\rm Ob}}
\newcommand \Aut {{\rm Aut}}
\newcommand \cor {{\rm {cor}}}
\newcommand \Cor {{\it {Corr}}}
\newcommand \res {{\rm {res}}}
\newcommand \red {{\rm{red}}}
\newcommand \Gal {{\rm {Gal}}}
\newcommand \PGL {{\rm {PGL}}}
\newcommand \Bl {{\rm {Bl}}}
\newcommand \Sing {{\rm {Sing}}}
\newcommand \spn {{\rm {span}}}
\newcommand \Nm {{\rm {Nm}}}
\newcommand \inv {{\rm {inv}}}
\newcommand \codim {{\rm {codim}}}
\newcommand \Div{{\rm{Div}}}
\newcommand \CH{{\rm{CH}}}
\newcommand \sg {{\Sigma }}
\newcommand \DM {{\sf DM}}
\newcommand \Gm {{{\mathbb G}_{\rm m}}}
\newcommand \tame {\rm {tame }}
\newcommand \znak {{\natural }}
\newcommand \lra {\longrightarrow}
\newcommand \hra {\hookrightarrow}
\newcommand \rra {\rightrightarrows}
\newcommand \ord {{\rm {ord}}}
\newcommand \Rat {{\mathscr Rat}}
\newcommand \rd {{\rm {red}}}
\newcommand \bSpec {{\bf {Spec}}}
\newcommand \Proj {{\rm {Proj}}}
\newcommand \pdiv {{\rm {div}}}
\newcommand \wt {\widetilde }
\newcommand \ac {\acute }
\newcommand \ch {\check }
\newcommand \ol {\overline }
\newcommand \Th {\Theta}
\newcommand \cAb {{\mathscr A\! }{\it b}}

\newenvironment{pf}{\par\noindent{\em Proof}.}{\hfill\framebox(6,6)
\par\medskip}

\newtheorem{theorem}[subsection]{Theorem}
\newtheorem{proposition}[subsection]{Proposition}
\newtheorem{lemma}[subsection]{Lemma}
\newtheorem{corollary}[subsection]{Corollary}

\theoremstyle{definition}
\newtheorem{definition}[subsection]{Definition}

\title[Divisibility of Selmer groups and class groups]{Divisibility of Selmer groups and class groups}
\author{Kalyan Banerjee, Kalyan Chakraborty, Azizul Hoque}
\address{Kalyan Banerjee @Harish-Chandra Research Institute, HBNI, Chhatnag Road, Jhunsi, Allahabad 211 019, India.}
\email{banerjeekalyan@hri.res.in}
\address{Kalyan Chakraborty @Harish-Chandra Research Institute, HBNI, Chhatnag Road, Jhunsi, Allahabad 211 019, India.}
\email{kalyan@hri.res.in}
\address{Azizul Hoque @Harish-Chandra Research Institute, HBNI, Chhatnag Road, Jhunsi, Allahabad 211 019, India.}
\email{ahoque.ms@gmail.com}


\begin{abstract}
In this paper, we study two topics. One is the divisibility problem of class groups of quadratic number fields and its connections to algebraic geometry. The other is the construction of Selmer group and Tate-Shafarevich group for an abelian variety defined over a number field.
\end{abstract}
\maketitle

\section{Introduction}
An interesting and important fact about the ideal class group of the ring of integers in a number field is that it is finite. So it is natural to ask whether given any positive integer $n$ there exists $n$-torsion elements in  the class group of the number field. This is related to the so called `divisibility problem' of the class group:  given a positive integer $n$ whether it divides the order of the class group. In the paper \cite{BH}, the authors studied the relations of this divisibility problem with the elements of $n$-torsions in the Picard group of a hyperelliptic surface. More precisely, suppose that we consider a hyperelliptic surface $S$ defined over $\bar {\QQ}$ given by a precise equation. Suppose that $S$ admits a regular map to the affine plane $\AF^2_{\bar {\QQ}}$ defined over $\bar \QQ$. Then spreading out $S,\AF^2$ over $\Spec(\ZZ)$ we have a family of ring of integers of a family of specific number fields. Now suppose that we start with a $n$-torsion element in the Picard group of $S$, then restricting it to fibers we have $n$-torsion elements in the class group of each member of the above family of ring of integers. Moreover a certain family of subgroups of $n$-torsions' of the class group of each member of the above family has same cardinality for each fiber over a Zariski open subgroup of $\AF^2_{\ZZ}$. This phenomena supports the so called Cohen-Lenstra heuristic which proposes the conjectural fact that given any $n$ a positive proportion of the number fields has an element  whcih are $n$-torsion in its class group. Atleast we can say that for the above family the number fields having the divisibility property is``parametrized" by a Zariski open subset of the affine plane over $\Spec(\ZZ)$.

The main theme of this work was to use the Mumford-Ro\u{i}tman argument for the natural map from relative Chow schemes to the relative Picard group, which says that the fibers of this map are countable union of Zariski closed subsets varying in a family.

The other theme of this paper is the Selmer and Tate-Shafarevich group constructions of an abelian variety defined over a number field.
The notion of Selmer group and the Tate-Shafarevich group is very much important from the perspective of local-global principle in arithmetic geometry. The Tate-Shafarevich group measures the failure of the local to global principle. The studies of these groups have been initiated by Cassels, Lang, Selmer, Shafarevich, Tate, \cite{CAS1}, \cite{CAS2}, \cite{LT}, \cite{Sel}, \cite{Sh1}, \cite{T}. The famous conjecture about the Tate-Shafarevich group tells that this group associated to an abelian variety is finite.  The first case where it has been proven is the case of elliptic curves with complex multiplication having rank atmost 1, by Karl Rubin, \cite{Ru1}. The next is the case of modular elliptic curves with analytic rank atmost 1, by V.Kolyvagin, \cite{Kol}. The paper by Selmer \cite{Sel} has many examples of genus one curves for which the Tate-Shafarevich group has non-trivial elements. The perception of the Tate-Shafarevich group of abelian varieties comes from the first Galois cohomology of the Abelian variety defined over a number field. On the other hand it can be described as the non-trivial torsors on the abelian variety which become trivial over a local field. The Selmer group has been defined by a certain kernel at the level of first Galois cohomology and it is known that this group is finite.

In the paper \cite{BC}, the  aim was to study the notion of Selmer group and the Tate-Shafarevich group from the perspective of algebraic cycles. That is the authors, consider the Galois action of the absolute Galois group of a number field on the group of degree zero cycles on an abelian variety defined over a number field. Then consider the Selmer and the Tate-Shafarevich group associated to this group of degree zero cycles on the abelian variety by considering the kernel at the level of first Galois cohomology of this particular Galois module.

There are certain things known about the group of degree zero cycles on a smooth projective variety over the algebraic closure of a number field. One is the Mumford-Ro\u{i}tman argument \cite{M},\cite{R1} about Chow schemes which says that the natural map from the symmetric power of an abelian variety to the Chow group has fibers given by a countable union of Zariski closed subsets in the symmetric power. This result enables us to give a scheme theoretic structure on the first Galois cohomology of the group of degree zero cycles on the abelian variety. Next, there is the famous theorem due to Ro\u{i}tman, \cite{R2}, which says that the torsion subgroup of the group of degree zero cycles on an abelian variety and the torsion subgroup on the abelian variety are isomorphic. This leads us to study of $n$-divisibility of the group of the degree zero cycles defined over the number field from a cohomological perspective. The main result of \cite{BC} is:
\begin{theorem}
Let $A$ denote an abelian variety defined over a number field $K$. Let $A_0(A)$ denote the group of degree zero cycles on the abelian variety. Then the group $A_0(A)(K)/nA_0(A)(K)$ is finite.
\end{theorem}
The importance of this result from the perspective of algebraic cycles lies in the Bloch-Beinlinson's conjecture on  the albanese kernel which says that the kernel of the albanese map from the group of degree zero cycles on a smooth projective variety over $\bar \QQ$ to the albanese variety has trivial kernel. It is known that this restriction on the ground field is sharp. That is, if we consider an one variable transcendental extension of the field $\bar \QQ$, then over this field there are varieties for which the albanese kernel is non-trivial (see \cite {GG},\cite{GGP}). From the above result we can only say that the quotient $T(A(K))/nT(A(K))$ of the Albanese kernel (denoted by $T(A)$) for an abelian variety $A$ is finite. Here $n$-is a positive integer and $A(K)$ denote the group of $K$-points on $A$. So it is worth studying the Tate-Shafarevich group and the Selmer group of the albanese kernel.

Here we have two constructions: one is the link between divisibility of class groups in a family with that of the Picard group of a hyperelliptic surface. On the other hand we can consider the Selmer group and the Tate-Shafarevich group of the Chow group of degree zero cycles on the hyperelliptic surface. A natural question is: What is the connection between $n$-divisibility of the Selmer group with that of the class group of the number fields in the above mentioned family?

The paper is organized as follows: in the second and third section we recall the results and techniques used in \cite{BC},\cite{BH} and then in the fourth section we attempt to connect these two techniques by answering the question posed in the last paragraph.

{\small \textbf{Acknowledgements:} K. Banerjee and A. Hoque thanks the hospitality of Harish Chandra Research Institute, India, for hosting this project. K.Banerjee was funded by DAE, Govt. of India, for this project and A. Hoque is supported by SERB N-PDF (PDF/2017/001958), Govt. of India.}

\section{Mumford-Ro\u{i}tman argument on Chow schemes and relative Picard schemes}
For a smooth, projective scheme $X$, let $D$ denote a Weil divisor  on it. For a flat morphism $X\to B$ of projective schemes, we consider the Chow scheme, $C^1_d(X/B)$ of relative divisors, that is co-dimension one subschemes of $X\to B$ of degree $d$, that is
$$C^1_{d}(X/B)=\{(D_b,b)|\Supp(D_b)\subset X_b, \deg(D_b)=d\}.$$
Then there is a natural map $C^1_d(X_b)\to \Pic(X_b)$ associating to a $D$ in its divisor class $[D]$. In this set up, we consider the following:
$$\bcZ:=\{(b,D_b)|[D_b]=0 \in \Pic(X_b)\}.$$
The proof of the following theorem is based on the idea introduced by Mumford in \cite{M}. This idea had been elaborated by Ro\u{i}tman in \cite{R} and Voisin in \cite{Voi}. This idea had also been used in \cite{BG}.

\begin{theorem}
\label{theorem1}
$\bcZ$ is a countable union of Zariski closed subsets in $C^1_{d}(X/B)$.
\end{theorem}

\begin{proof}
Assume that the relation $D_b=D_b^+-D_b^-$ is rationally equivalent to zero. This means that there exists a map $f:\PR^1\to C^1_{d,d}(X_b)$ such that
$$f(0)=D_b^{+}+\gamma\text{ and }f(\infty)=D_b^{-}+\gamma,$$
where $\gamma$ is a positive divisor on $X_b$.
In other words, we have the following map: $$\ev:Hom^v(\PR^1_k,C^1_{d}(X/B))\to C^1_d(X/B)\times C^1_d(X/B),$$ given by $f\mapsto (f(0),f(\infty))$ and image of $f$ is contained in $C^1_{d,d}(X_b)$.

Let us denote $C^1_d(X/B)$ by $C^1_d(X)$ for simplicity.

We now consider the subscheme $U_{v,d}(X)$ of $B\times Hom^v(\PR^1_k,C^1_{d}(X))$ consisting of the pairs $(b,f)$ such that image of $f$ is contained in $C^1_{d}(X_b)$ (such a universal family exists, for example see \cite[Theorem 1.4]{Kol}). This gives a morphism from $U_{v,d}(X)$ to $B\times C^1_{d,d}(X)$ defined by $$(b,f)\mapsto (b,f(0),f(\infty)).$$
Again, we consider the closed subscheme $\bcV_{d,d}$ of $B\times C^1_{d,d}(X)$ given by $(b,z_1,z_2)$, where $(z_1,z_2)\in C^1_{d,d}(X_b)$. Suppose that the map from $\bcV_{d,u,d,u}$ to $\bcV_{d+u,d+u}$ is given by
$$(A,C,B,D)\mapsto (A+C,C,B+D,D).$$
Then one writes the fiber product $\bcV$ of $U_{v,d}(X)$ and $\bcV_{d,u,d,u}$ over $\bcV_{d+u,d+u}$. If we consider the projection from $\bcV$ to $B\times C^1_{d,d}(X)$, then we observe that $A$ and $B$ are supported as well as rationally equivalent on $X_b$. Conversely, if $A$ and $B$ are supported as well as rationally equivalent on $X_b$, then one gets the map $$f:\PR^1\to C^1_{d+u,u,d+u,u}(X_b)$$ of some degree $v$ satisfying
$$f(0)=(A+C,C)\text{ and } f(\infty)=(B+D,D),$$
where $C$ and $D$ are supported on $X_b$. This implies that the image of the projection from $\bcV$ to $B\times C^1_{d,d}(X)$ is a quasi-projective subscheme $W_{d}^{u,v}$ consisting of the tuples $(b,A,B)$ such that $A$ and $B$ are supported on $X_b$, and that there exists a map $$f:\PR^1_k\to C^1_{d+u,u}(X_b)$$ such that $f(0)=(A+C,C)$ and $f(\infty)=(B+D,D)$. Here $f$ is of degree $v$, and $C,D$ are supported on $X_b$ and they are of co-dimension $1$ and degree $u$ cycles. This shows that $W_d=\cup_{u,v} W_d^{u,v}$. We now prove that the Zariski closure of $W_d^{u,v}$ is in $W_d$ for each $u$ and $v$. For this, we prove the following:
$$W_d^{u,v}=pr_{1,2}(\wt{s}^{-1}(W^{0,v}_{d+u}\times W^{0,v}_u)),$$
where
$$\wt{s}: B\times C^1_{d,d,u,u}(X)\to B\times C^1_{d+u,d+u,u,u}(X)$$
defined by
$$\wt{s}(b,A,B,C,D)=(b,A+C,B+D,C,D).$$

We assume $(b,A,B,C,D)\in B\times C^1_{d,d,u,u}(X)$ in such a way that $\wt{s}(b,A,B,C,D)\in W^{0,v}_{d+u}\times W^{0,v}_u$. This implies that there exists an element $(b,g)\in B\times\Hom^v(\PR^1_k,C^p_{d+u}(X))$ and an element $(b,h)\in \Hom^v(\PR^1_k,C^p_{u}(X))$ satisfying $$g(0)=A+C,~g(\infty)=B+D \text{ and } h(0)=C,h(\infty)=D$$ as well as the image of $g$ and $h$ are contained in $C^1_{d+u}(X_b)$ and  $C^1_u(X_b)$ respectively.

Also if $f=g\times h$ then $f\in \Hom^v(\PR^1_k,C^p_{d+u,u}(X))$ such that the image of $f$ is contained in $C^1_{d+u,u}(X_b)$ as well as it satisfies the following:
$$f(0)=(A+C,C)\text{ and }(f(\infty))=(B+D,D).$$
This shows that $(b,A,B)\in W^d_{u,v}$.

On the other hand, we assume that $(b,A,B)\in W^d_{u,v}$. Then there exists $f\in \Hom^v(\PR^1_k,C^1_{d+u,u}(X_b))$ such that
$$f(0)=(A+C,C)\text{ and }f(\infty)=(B+D,D),$$
and image of $f$ is contained in the Chow scheme of  $X_b$.

We now compose $f$ with the projections to $C^1_{d+u}(X_b)$ and to $C^1_{u}(X_b)$ to get a map $g\in \Hom^v(\PR^1_k,C^1_{d+u}(X))$ and a map $h\in\Hom^v(\PR^1_k,C^1_{u}(X))$ satisfying
$$g(0)=A+C,~g(\infty)=B+D$$
and
$$h(0)=C,~h(\infty)=D.$$
Also, the image of $g$ and $h$ are contained in the respective Chow schemes of   the fibers $X_b$. Therefore, we have
$$W_d=pr_{1,2}(\wt{s}^{-1}(W_{d+u}\times W_u)).$$
We are now in a position to prove that the closure of $W_d^{0,v} $ is contained in $W_d$. Let $(b,A,B)$ be a closed point in the closure of ${W_d^{0,v}}$. Let $W$ be an irreducible component of ${W_d^{0,v}}$ whose closure contains $(b,A,B)$. Assume that $U$ is an affine neighborhood of $(b,A,B)$ such that $U\cap W$ is non-empty. Then there is an irreducible curve $C$ in $U$ passing through $(b,A,B)$. Suppose that $\bar{C}$ is the Zariski closure of $C$ in $\overline{W}$. The map
$$e:U_{v,d}(X)\subset B\times \Hom^v(\PR^1_k,C^1_{d}(X))\to C^1_{d,d}(X)$$
given by
$$(b,f)\mapsto (b,f(0),f(\infty))$$
is regular and $W_d^{0,v}$ is its image. We now choose a curve $T$ in $U_{v,d}(X)$ such that the closure of $e(T)$ is $\bar C$.  Let $\wt{T}$ be denote the normalization of the Zariski closure of $T$, and $\wt{T_0}$ be the pre-image of $T$ in this normalization. Then the regular morphism $\wt{T_0}\to T\to \bar C$ extends to a regular morphism from $\wt{T}$ to $\bar C$. If $(b,f)$ is a pre-image of $(b,A,B)$, then $f(0)=A,~ f(\infty)=B$ and the image of $f$ is contained in $C^p_{d}(X_b)$ by the definition of $U_{v,d}(X)$. Therefore, $A$ and $B$ are  rationally equivalent. This completes the proof.
\end{proof}

As a consequence, one gets the following:

\begin{corollary}\label{cor2.1} The collection
$$\bcZ_d:=\{(b,D_b)|n[D_b]=0 \in \Pic(X_b)\}$$
is a countable union of Zariski closed subsets in the scheme $C^1_d(X/B)$.
\end{corollary}

\subsection{Mumford-Ro\u{i}tman arguments and monodromy representation}

Following an idea using monodromy due  to Voisin \cite[Chapter 3]{Vo} and the above mentioned argument due to Mumford and Ro\u{i}tman, we have the following theorem:

\begin{theorem} \label{theorem2}
The cardinality of the subgroup of torsions in $X_b$ coming from the fibration $\bcZ_{i\CC,U}\to U_{\CC}$ for each $b\in U$ remains constant and they vary in a family.
\end{theorem}
\begin{proof}
For a proof see \cite{BH}[Thoerem 3.1].
\end{proof}

We now consider a smooth projective curve $C$ over an algebraically closed field $K\subset \CC$ in the  projective plane $\PR^2$ over $K$. Let $U$ be an affine piece of $C$. That is, $U$ is $C$ minus finitely many points, viz. $P_1,\cdots,P_m$. Consider the following localization exact sequence of Picard groups
$$\oplus_i \ZZ[P_i]\to \Pic(C)\to \Pic(U)\to 0.$$
Then the set of all torsion points in $\Pic(U)$ gives rise to elements of $\Pic(C)$ of the form $nz$ such that
$$nz=\sum_i n_iP_i$$
where $P_1,\cdots, P_m$ are the finite number of points that are deleted. As before, we consider a fibration of smooth projective schemes $X\to B$ over $\bar \QQ$, where $X$ is a surface embedded in $\PR^3$ such that each fiber $X_b$ is contained in a projective plane $\PR^2$ over $\bar \QQ$ and $B$ is an algebraic curve. Suppose that the degree of the algebraic curve $X_b$ remains constant over a Zariski open set $U$ in $B$. For an affine piece $U_b$ of the algebraic curve $X_b$, we consider the following:
$$\bcP:=\{(x,b)|x\in X_b\setminus U_b\}\to U.$$
By the above assumption, this a finite-to-one map from $\bcP$ to $U$ and the degree of this map is constant. For  given any $b\in U$, let us suppose the fiber $\bcP_b$ contains the points $P_{1b},\cdots,P_{mb}$. We define the set:
$$\bcZ_d=\{(b,z)|\Supp(z)\subset X_b, nz=\sum_i n_iP_{ib}\}\to U.$$
Then as a conseqeunce of Theorem \ref{theorem1} one gets the following result.
\begin{corollary}
$\bcZ_d$ is a countable union of Zariski closed subsets in the ambient relative Chow scheme $C^1_d(X_U/U)$, where $X_U\to U$ is the pullback of the family $X\to B$ to $U$.
\end{corollary}
This corollary along with the Theorem \ref{theorem2} gives the following:
\begin{corollary} \label{theorem3}
The cardinality of the set of $z$ in $\bcZ_{ib}$ for $b\in U$ such that
$$nz=\sum_i n_i P_{ib}$$
for points $P_{ib}\in X_b$ is constant as $b$ varies over $U$.
\end{corollary}
These points on $\Pic(X_b)$ correspond to the torsion elements in $\Pic(U_b)$, where $U_b$ is the open complement of $X_b$ obtained from $X_b$ by  deleting the points $P_{1b},\cdots, P_{mb}$.

\subsection{Example of hyperelliptic surfaces}

In this section, we will show that certain algebraic surfaces have $n$-torsion elements in the Picard group. We begin section with the algebraic surface defined by
$$y^2=t^2q^2-z^n$$
over $\QQ$. Its co-ordinate ring is given by
$$\QQ[y,t,z]/(y^2-t^2q^2+z^n).$$
We now consider the maximal ideal $(t-m,z-\ell)$, for some  algebraic  numbers $m,\ell$,  in the polynomial ring $\QQ[t,z]$. We also consider the map $$\QQ[t,z]\to \QQ[y,t,z]/(y^2-t^2q^2+z^n)$$
which is defined by
$$t\mapsto t, z\mapsto z$$
and the map $\QQ[t,z]\to \QQ$ which is given by
$$f(t,z)\mapsto f(m,\ell).$$
Then the tensor product
$$ \QQ[y,t,z]/(y^2-t^2q^2+z^n)\otimes _{\QQ[t,z]}\QQ$$
is given by $ \QQ[y]/(y^2-m^2q^2+\ell^n)$. Further, if the polynomial $p(y):=y^2-m^2q^2+\ell^n$ is irreducible over $\QQ$, then the above co-ordinate ring is isomorphic to $L$, where $L$ is the imaginary quadratic extension of $\QQ$ given by adjoining a root of $p(y)$. Therefore if we consider the family $$\ZZ[y,t,z]/(y^2-t^2q^2+z^n)\to \ZZ[t,z],$$
then the normalizations of the fibers are the ring of integers of 
$$\QQ(\sqrt {m^2q^2-\ell^n})$$

Let us consider an affine surface $S$ fibered over $\AA^2_{\QQ}$ as mentioned in the beginning of this section. Let the pullback of the fibration over $\AA^2_{\ZZ}$, be $S_{\ZZ}\to \AA^2_{\ZZ}$. Then the family of mormalizations  to $S_{\ZZ}$ is the family of ring of integers $\bcO(\sqrt{m^2q^2-\ell^n})$. For the convenience of notation, let us continues to denote this family as $S_{\ZZ}$. Consider the Zariski closure of $S$ in $\PR^3_{\QQ}$ and the Zariski closure of the family $S_{\ZZ}\to \AA^2_{\ZZ}$ in $\PR^3_{\ZZ}$. We denote it by $\bar S_{\ZZ}$. We also consider the Chow scheme
$$C^1_d(\bar S_{\ZZ}/\PR^2_{\ZZ})$$
and the subset
$$\bcZ_d:=\{(z,b)|\Supp(z)\subset X_b, [z]=\sum_i n_i[P_{ib}]\in \Pic(S_{\ZZ,b})\},$$
where $P_{1b},\cdots,P_{mb}$ are the points in the complement of $S_{\ZZ,b}$ inside the Zariski closure $\bar S_{\ZZ,b}$. Then by Theorem \ref{theorem1}, we get the following result.
\begin{proposition}
The set $\bcZ_d$ is a countable union of Zariski closed subsets in the Chow scheme.
\end{proposition}
Applying the same argument as in Corollary \ref{theorem3}, we see that there exists an irreducible Zariski closed subset $\bcZ_i$ inside the relative Picard scheme $\Pic(\overline{S_{\ZZ U}}\to U)$, where $U$ is Zariski open in $\AA^2_{\ZZ}$, such that the complexification of $\bcZ_{i,\CC}$ maps dominantly onto $U_{\CC}$ as well as the number of points in the fiber of this map is constant. Therefore, one gets the following:
\begin{theorem}\label{thm3.5}
The cardinality of a certain subgroup of $\Pic(S_{\ZZ,b})$ which is nothing but the class group of the quadratic field $\QQ(\sqrt{m^2q^2-\ell^n})$ for some fixed integers $m$ and $\ell$, remains constant as $b$ varies over $U$.
\end{theorem}
This concludes that given an element of order $n$ in $\Pic(S_{\ZZ,b})$, one can find an element of the same order in $\Pic(S_{\ZZ,b'})$ for some $b'\in U$ which is different from $b$.

\section{Tate-Shafarevich group of the Chow group of an abelian variety}

Let $K$ be a number field and let $\overline{K}$ denote its algebraic closure. Let $A$ be an abelian variety defined over  $K$. Then we have a natural Galois action of the absolute Galois group $\Gal(\bar K/K)$. This action induces further an action on the Chow group of zero cycles on the abelian variety $A$. Here the Chow group is the free abelian group generated by closed points on $A(\bar K)$ modulo the rational equivalence. We denote this group by $\CH_0(A)$.

Consider the continuous functions $f$ from $G=\Gal(\bar K/K)$ to $\CH_0(A)$ satisfying the property that
$$f(\sigma\tau)=f(\sigma)+\sigma f(\tau)\;.$$
The set of all such functions form a group denoted by $Z^1(G,\CH_0(A))$. Let us consider the subgroup of $Z^1(G,\CH_0(A))$ consisting of elements $f$ such that
$$f(\sigma)=\sigma.x-x$$
where $x$ some element in the group $\CH_0(A)$. Denote this subgroup by $B^1(G,\CH_0(A))$. Then we define the quotient
$$Z^1(G,\CH_0(A))/B^1(G,\CH_0(A))$$
as
$$H^1(G,\CH_0(A))\;.$$

Let $A_0(A)$ denote the subgroup of degree zero cycles modulo rational equivalence in $\CH_0(A)$. We consider as previous the Galois cohomology $$H^1(G,A_0(A))$$
of the group $A_0(A)$.

We observe that there is a natural homomorphism of abelian groups from $A_0(A)$ to $A$. Then by functoriality of group cohomology we have that this homomorphism descends to a homomorphism of Galois cohomology groups of the corresponding Galois modules: 
$$H^1(G,A_0(A))\to H^1(G,A)$$ 
The map from $A_0(A)$ to $A$ is denoted by $\alb$, the albanese map. We denote the map from $H^1(G,A_0(A))$ to $H^1(G,A)$ as $\alb$. We are interested in understanding the structure of the group $H^1(G,A_0(A))$. Consider the natural map from $\Sym^n A$ to $A_0(A)$, which sends an unordered $n$-tuple $\{P_1,\cdots,P_n\}$ of $\bar K$ points on $A$ to the cycle class
$$\sum_{i=1}^n [P_i-P_0]\;,$$
where $P_0$ is a fixed $K$-point on $A$.

Now consider the fact that the group $H^1(G,A_0(A))$ is actually isomorphic to the colimit of Galois cohomology of finite groups $$H^1(\Gal(L/K),A_0(A_L))$$ Here $L/K$ is a finite Galois extension and $A_L$ is the collection of $L$-points on $A$. Since $G$ is a profinite group, the range of any function $\eta$ from $G$ to $A_0(A)$ is finite.  Consider $Z_{l}$ to be the collection of all maps $\eta$ from $G$ to $A_0(A)$ such that $\eta$ factors through $\Sym^l A\times \Sym^l A$ (this can be achieved by decomposing a zero cycle into positive and negative parts). That is we identify the maps $\eta$, factoring through $\Sym^l A\times \Sym^l A$, with its image inside $\Sym^l A\times \Sym^l A$. There exists a normal subgroup of $G$ of finite index, call it $N$, such that $\eta$ is factoring through $G/N$. On the other hand suppose that we have a collection of points on $\Sym^l A\times \Sym^l A$. Then we can define a map from $G/N$ to $\Sym^l A\times \Sym^l A$ by assigning the cosets of $N$ to this finite collection of points of $\Sym^l A\times \Sym^l A$. Such a map will be continuous from $G/N$ to $\Sym^l A\times \Sym^l A$ equipped with discrete topology, as $G/N$ is finite. Since the quotient map from $G$ To $G/N$ is continuous we have that the map from $G$ to $\Sym^l A\times \Sym^l A$ is continuous. But these maps are non-canonical as its depends on the choice of the points and their assignments to the left cosets of $N$.  Now consider the relation that defines $Z^1(G,A_0(A))$,
$$\eta(\sigma \tau)=\eta(\sigma)+\sigma\eta(\tau)\;.$$
Since this relation happens on $A_0(A)$ we have that the cycles
$$\eta(\sigma\tau)$$
is rationally equivalent to
$$\eta(\sigma)+ \sigma\eta(\tau)\;.$$
This means that there exists a map from $\PR^1_{\bar K}$ to $\Sym^d A$, and a positive zero cycle $B$ such that
$$f(0)=\eta(\sigma\tau)+B, \quad f(\infty)=\eta(\sigma)+\sigma\eta(\tau)+B\;.$$
By the theorem of Roitman \cite{R} the collection of all such $\eta$, such that
$$\eta(\sigma\tau)$$
{is rationally equivalent to} $$\eta(\sigma)+\sigma\eta(\tau)$$
is a countable union of Zariski closed subsets inside the symmetric power $\Sym^l A\times \Sym^l A$ such that range of $\eta$ is contained in $\Sym^l A\times \Sym^l A$.  So following \cite{R} we have:
\begin{theorem}
The collection of all $\eta$ contained in $Z_{l,l}$ such that
$$\eta(\sigma\tau)$$
is rationally equivalent to
$$\eta(\sigma)+\sigma\eta(\tau)$$
is a countable union of Zariski closed subsets inside $\Sym^l A\times \Sym^l A$ denoted by $Z^1_{l}$.
\end{theorem}
\begin{proof}
For details of the proof  see \cite{BC}[theorem 2.1].
\end{proof}

Similarly we can prove that the collection of $\eta$ in $Z_{l}$ such that $\eta(\sigma)$ is rationally equivalent to $\sigma.z-z$ (for a fixed zero cycle $z$) is a countable union of Zariski closed subsets in $Z^1_{l}$. Call it $B^1_{l}$.

Therefore we can conclude from the above theorem that:

\begin{theorem}
The group $H^1(G,A_0(A))$ admits a surjective map from the countable union $\cup_{l}Z^1_{l}$ such that $\cup_{l}B^1_{l}$ is mapped to a point under this surjective map.
\end{theorem}

Now we further study the property of this map: 
$$Z^1_{l}\to H^1(G,A_0(A))$$ 
Consider an element $\eta$ in the set $Z^1_{l}$. Then for every $\sigma,\tau$ we have
$$\eta(\sigma\tau)=\eta(\sigma)+\sigma\eta(\tau)\;.$$
This equality happens in $A_0(A)$. So consider the tuples
$$(\eta,f,B)\in Z_{l}\times \Hom^v(\PR^1,\Sym^{n+u,n+u}A)\times \Sym^u B$$
such that the following equations are satisfied:
$$f(0)=\eta(\sigma\tau)+B$$
$$f(\infty)=\eta(\sigma)+\sigma\eta(\tau)+B.$$
So if we denote the above quasiprojective variety by $\bcV$ and consider the projection map from $\bcV$ to $\Hom^v(\PR^1,\Sym^{n+u,n+u}A)$, then it is a $\PR^1$-bundle. This is because it is the pull-back of the $\PR^1$-bundle given by
$$\{(x,f)|x\in \im (f)\}\subset \Sym^{n+u,n+u}A\times \Hom^v(\PR^1,\Sym^{n+u,n+u}A)\;.$$
So over $Z^1_{l}$ we have the universal variety $\bcU^1_{l,m}$ consisting of tuples $(\eta,f,B)$ such that the above equations are satisfied and it has the structure of a rationally connected fibration over the Hom-scheme. Therefore if we consider the finite map from $A^{2l}$ to $\Sym^l A\times \Sym^l A$, the degree of this finite map is $(l!)^2$. The pullback of $Z^1_{l}$ under this map is a finite branched cover of $Z^1_{l}$ denoted by $\widetilde{Z^{1}_{l}}$. Correspodningly we have the pull-back of the universal family $\bcU^1_{l}$ over $\widetilde{Z^{1}_{l}}$ denoted by $\widetilde{\bcU^{1}_{l}}$. This is a family of branched covers of $\PR^1$ over the Hom-scheme.

\subsection{The group cohomology of the group of degree zero cycles on A}
Let $A_0(A)$ denote the group of degree zero cycles or the zero cycles algebraically equivalent to zero on $A$. Then there is a natural homomorphism from $A^n$ to $A_0(A)$ Given by
$$\sum_i P_i\mapsto \sum_i [P_i-n0]$$
here $0$ is the neutral element of the abelian variety $A$. Then the map from $A^n$ to $A_0(A)$ induces by functoriality a natural homomorphism from $H^1(G,A^n)$ to $H^1(G,A_0(A))$. Consider the natural map from $A^n$ to $A^{n+1}$ given by
$$(P_1,\cdots,P_n)\mapsto (P_1,\cdots,P_n,0)$$
Then this map gives rise to the homomorphism from $H^1(G,A^n)$ to $H^1(G,A^{n+1})$ and the homomorphism
$$\theta_n:H^1(G,A^n)\to H^1(G,A_0(A))$$
factors through the above map
$$H^1(G,A^n)\mapsto H^1(G,A^{n+1})\;.$$
Hence we have a natural homomorphism from the colimit of the groups
$$H^1(G,A^n)$$ to
$$H^1(G,A_0(A))$$ denoted by $\theta$. So we have
$$\theta:\varinjlim H^1(G,A^n)\to H^1(G,A_0(A))\;.$$
Now for each $n$ we have the group law from $A^n$ to $A$ given by
$$(a_1,\cdots,a_n)\mapsto \sum_i a_i$$
This map gives rise to a natural map from $H^1(G,A^n)$ to $H^1(G,A)$. Note that this map factors through the homomorphism
$$H^1(G,A^n)\to H^1(G,A^{n+1})$$
Therefore we have a homomorphism from
$$\varinjlim H^1(G,A^n)\to H^1(G,A)\;.$$
Since the map $H^1(G,A^n)$ to $H^1(G,A_0(A))$ factors through the map
$$H^1(G,A)\to H^1(G,A_0(A))$$
we have that the map
$$\varinjlim H^1(G,A^n)\to H^1(G,A_0(A))$$
factors through the map
$$H^1(G,A)\to H^1(G,A_0(A))\;.$$
Now the group on the left is the the Weil-Chatelet group of the respective $A$, which consists of the equivalence classes of principal homogeneous spaces over $A$. This group is denoted by $WC(A)$. Under the identification
$$H^1(G,A)\cong WC(A)$$
we have that
$$\theta: WC(A)\to H^1(G,A_0(A))\;.$$
It is natural to consider when this map is injective and surjective.

Now due to the famous result on torsions by Ro\u{i}tman \cite{R2} in $A_0(A)$, we know that this group of torsions is isomorphic to the group of torsions in $A$. So we expect a similar result when we consider the group cohomology
$H^1(G,A)$ and $H^1(G,A_0(A))$.

\begin{theorem}
The kernel of the map $H^1(G,A)[n]\to H^1(G,A_0(A))[n]$ is isomorphic to the group
$$A_0(A(K))/nA_0(A(K))\;.$$
\end{theorem}

\begin{proof}
For a proof see \cite{BC}[theorem 2.4].
\end{proof}

\subsection{Tate-Shafarevich and Selmer group of $A_0(A)$ and their properties}

Consider the exact sequence

$$0\to A_0(A(K))/nA_0(A(K))\to H^1(G,A_0(A)[n])\to H^1(G,A_0(A))[n]\to 0\;.$$

Now we consider a place $v$ of $K$ and consider the complection of $K$ at $v$, denote this completion by $K_v$. Then consider the algebraic closure $\bar K_v$ of $K_v$ and embed $\bar K$ into $\bar K_v$. This embedding gives us an injection of the Galois group $\Gal(\bar K_v/K_v)$ into $\Gal(\bar K/K)$. Considering the Galois cohomology, we have a homomorphism from
$$H^1(\Gal(\bar K/K),A_0(A(\bar K))\to H^1(\Gal(\bar K_v/K_v),A_0(A(\bar K_v)))\;.$$
We write the groups $\Gal(\bar K/K),\Gal(K_v/K_v)$ as $G,G_v$ for simplicity.
Then we have the following commutative diagrams:

$$
  \diagram
  A_0(A(K))/nA_0(A(K))\ar[dd]_-{} \ar[rr]^-{} & & H^1(G,A_0(A)[n]) \ar[dd]^-{} \\ \\
  A_0(A(K_v))/nA_0(A(K_v)) \ar[rr]^-{} & & H^1(G_v,A_0(A_v)[n])
  \enddiagram
  $$

$$
  \diagram
  H^1(G,A_0(A)[n])\ar[dd]_-{} \ar[rr]^-{} & & H^1(G,A_0(A))[n] \ar[dd]^-{} \\ \\
 H^1(G_v,A_0(A_v)[n])\ar[rr]^-{} & & H^1(G_v,A_0(A_v))[n]
  \enddiagram
  $$
Consider the map
$$H^1(G,A_0(A)[n])\to \prod_v H^1(G_v,A_0(A_v))$$

\begin{definition}
The kernel of this map is defined to be the selmer group associated to the map $z\mapsto nz$ and it is denoted by $S^n(A_0(A)/K)$.
\end{definition}

\begin{definition}
The Tate-Shafarevich group is the kernel of the map
$$H^1(G,A_0(A))\to \prod_v H^1(G_v,A_0(A_v))\;.$$
It is denoted by $TS(A_0(A)/K)$.
\end{definition}

Now consider the commutative diagram:

$$
  \diagram
  H^1(G,A[n])\ar[dd]_-{} \ar[rr]^-{} & & \prod_v H^1(G_v,A_v)[n] \ar[dd]^-{} \\ \\
 H^1(G,A_0(A)[n])\ar[rr]^-{} & & \prod_v H^1(G_v,A_0(A_v))[n]
  \enddiagram
  $$
Now by Ro\u{i}tman's theorem as in \cite{R2}, the groups $A[n]$ and $A_0(A)[n]$ are isomorphic, as Galois modules (after a possible finite extension of the given number field). This fact is explained in details in the next section \ref{section4}. Therefore the group cohomologies are isomorphic. So the left vertical arrow in the above diagram is an isomorphism. Suppose that we have an element in $S^n(A/K)$, then by the commutativity of the above diagram we have that the image of the element under the left vertical homomorphism is in $S^n(A_0(A)/K)$. Now we recall the following theorem proved in \cite{BC}[theorem 3.3]

\begin{theorem}
The group $S^n(A_0(A)/K)$ is finite and hence
$$A_0(A(K))/nA_0(A(K))$$
 is finite.
\end{theorem}

\section{Divisibility problem of class groups and Selmer groups}
\label{section4}

Considering the group of $n$-torsion elements in the Picard variety  of the hyperelliptic surface $S$, say $\Pic^0(S)$, we have a bijection of this group with the $n$-torsion subgroup of $\Alb(S)$, the albanese variety of $S$. This isomrophism is defined over a finite extension of the ground field. On the other hand by the Ro\u{i}tman's theorem  in \cite{R2}, the $n$-torsions in $A_0(S)$ corresponds to $n$-torsions' in $\Alb(S)$. Suppose that after a finite extension of the ground field, say $K$, we have a $K$-rational point on $S$. Call it $P_0$. Let $g$ be an element in the absolute Galois group of $K$. Then we have a functorial morphism:
$$g_*:\Sym^n S \to \Sym^n S$$
given by 
$$P_1+\cdots+P_n\mapsto g(P_1)+\cdots+g(P_n)$$
where $P_1+\cdots+P_n$ denote the unordered $n$-tuple consisting of closed points $P_1,\cdots,P_n$ in $S$.
Consider the natural map 
$$\Sym^n S\to A_0(S)$$
given by
$$P_1+\cdots+P_n\mapsto [P_1+\cdots+P_n-nP_0]$$
The right hand side above denote the cycles class corresponding to the cycle
$$\sum_{i=1}^n P_i-nP_0\;.$$
Note that the above map gives rise to the formula:
$$g_*([P_1+\cdots+P_n-nP_0])=[g_*(P_1)+\cdots+g_*(P_n)-ng_*(P_0)]$$
for $g$ an element in the Galois group. Therefore we have 
$$g_*:A_0(S)\to A_0(S)$$
composing this map with the albanese map we have 
$$A_0(S)\to A_0(S)\to \Alb(S)\;.$$
If we compose the map $S\to A_0(S)$, with the above map then we have 
$P_0$ mapping to zero in $\Alb(S)$. Hence we have a unique morphism of abelian varieties (denote it again by $g_*$):
$$\Alb(S)\to \Alb(S)$$
such that we have 
$$\alb_S\circ g_*=g_*\circ \alb_S$$
Here $\alb_S$ is the albanese map. This gives that the map 
$$A_0(S)\to \Alb(S)$$
is a map of Galois modules (ensured by the universal property of the albanese variety), provided that $S$ has a $K$-rational point.
Let for an abelian group $A$, the $n$-torsions are denoted by $A[n]$. So we have the isomorphisms of Galois modules
$$\Pic^0(S)[n]\to \Alb(S)[n]\to A_0(S)[n]$$
where the fist one comes from Autoduality of Picard and Albanese varieties and the second one is the isomorphism coming from Ro\u{i}tman's theorem described above. Then we have an isomorphism between
$$H^1(G,\Pic^0(S)[n])$$
and
$$H^1(G,A_0(S)[n])\;.$$

Therefore considering the Tate-Shafarevich groups we have an isomorphism
$$TS(A_0(S)[n]/K)\cong TS(\Pic^0(S)[n]/K)$$
where $TS$ denote the Tate-Shafarevich group of $n$-torsions of the corresponding group. Now suppose that we start from an element of
order $n$ on $A_0(S)$, this will correspond to an element of order $n$ in the Selmer group in $\Pic^0(S)$. Now consider a fibration of this surface over $\PR^1_{\bar \QQ}$. Then for a closed point $b\in \PR^1_{\bar \QQ}$, we have the restriction homomorphism
$$\Pic^0(S)\to \Pic^0(S_b)$$

So considering both $\Pic^0(S),\Pic^0(S_b)$ as Galois modules we have a map of Galois cohomology
$$H^1(G,\Pic^0(S))\to H^1(G,\Pic^0(S_b))$$
Thus we have the following commutative diagram:

$$
  \diagram
 \Pic^0(S)(K)/n\Pic^0(S)(K)\ar[dd]_-{} \ar[rr]^-{} & & H^1(G,\Pic^0(S)[n])\ar[dd]^-{} \\ \\
\Pic^0(S_b)(K)/n\Pic^0(S_b)(K)\ar[rr]^-{} & & H^1(G,\Pic^0(S_b)[n])
  \enddiagram
  $$
Also we have an analogue of the above diagram at the level of local fields that is $K_v$, where $v$ is a place of $K$. Therefore functorially we have the homomorphism from
$$TS(\Pic^0(S)/K)\to TS(\Pic^0(S_b)/K)$$
and
$$S^n(\Pic^0(S)/K)\to S^n(\Pic^0(S_b)/K).$$
Composing this map with the map
$$S^n(A_0(S)/K)\to S^n(\Pic^0(S)/K)$$
we have established a homomorphism from $S^n(A_0(S)/K)$ to $S^n(\Pic^0(S_b)/K)$ and corresponding homomorphism
$$TS(A_0(S)/K)\to TS(\Pic^0(S_b)/K)\;.$$
By a diagram chase we have the following diagrams
$$
  \diagram
 A_0(S)(K)/nA_0(S)(K)\ar[dd]_-{} \ar[rr]^-{} & & S^n(A_0(S)/K)\ar[dd]^-{} \\ \\
\Pic^0(S_b)(K)/n\Pic^0(S_b)(K)\ar[rr]^-{} & & S^n(\Pic^0(S_b))
  \enddiagram
  $$

\smallskip

$$
  \diagram
 S^n(A_0(S)/K)\ar[dd]_-{} \ar[rr]^-{} & & TS^n(A_0(S)/K)\ar[dd]^-{} \\ \\
S^n(\Pic^0(S_b))\ar[rr]^-{}  & & TS^n(\Pic^0(S_b)/K)
  \enddiagram
  $$
Therefore the elements of order $n$ in the Selmer or the Tate-Shafarevich group of $A_0(S)$ corresponds to an element of order $n$ in the Selmer or the Tate-Shfarevich group respectively of $\Pic^0(S_b)$. So the $n$-divisibility of the Selmer or the Tate-Shafarevich group of $A_0(S)$ corresponds to $n$-divisibility of the Selmer or the Tate-Shafarevich group of $\Pic^0(S_b)$. Now spreading out $S_b$ over $\Spec(\ZZ)$ and considering an affine piece of the spread we have the $\Spec(\bcO_b)$, where $\bcO_b$ is the ring of integers of a number field. The above construction then corresponds to the following:
\begin{theorem}
The $n$-divisibility of the Tate-Shafarevich or the Selmer group of the group $A_0(S_U)$ of an affine piece $S_U$ of the surface $S$, to the respective $n$-divisibility of the corresponding class group of $\bcO_b$.
\end{theorem}


\begin{thebibliography}{AAAAA}
\bibitem[BC]{BC} K. Banerjee and K. Chakraborty, {\em Tate-Shafarevich group and Selmer group constructions for Chow group of an abelian variety },
{\small \tt arxiv: 1906.08233}, 2019.
\bibitem[BG]{BG} K. Banerjee and V. Guletskii, {\em \'{E}tale monodromy and rational equivalence for 1-cycles on cubic hypersurfaces in $\mathbb P^5$}, Mat. Sb. (to appear) {\small \tt arXiv:1405.6430v1}.
\bibitem[BH]{BH} K. Banerjee and A. Hoque, {\em Picard group, pull back and class group}, {\small \tt arxiv: 1903.04210}, 2019.
\bibitem[CAS1]{CAS1} J. Cassels, {\em Arithmetic of curves of genus 1. III. Tate-Shafarevich and Selmer groups }, Proc. London Math. Soc. Third series, 12, 259--296
\bibitem[CAS2]{CAS2} J. Cassels, {\em Arithmetic of curves of genus 1. IV. Proof of Hauptvermutung}, J. Reine Angew. Math. 211, 95--112.
\bibitem[GG]{GG} S. Gorchinsky and V. Guletskii, {Non-trivial elements in the kernel of Abel-Jacobi kernels of higher dimensional varities}, Adv. Math. 241, 2013, 162--191.
\bibitem[GGP]{GGP} M. Green, P. Griffiths and K. Paranjape, {Cycles over fields of transcendence degree one}, Michigan Math. J. 52(1), 2004.
\bibitem[Kol]{Kol} V.Kolyvagin, {\em Finiteness of E(Q) and SH(E(Q)) for a subclass of Weil-Curves}, Izvestiya Akademii Nauk SSSR, Seriya Matematicheskaya, 52(3), 522-540, 670-671, 1988.
\bibitem[LT]{LT} S.Lang, J.Tate, {\em Principal homogeneous spaces over abelian varities}, Americal journal of Mathematics, 80(3), 659-684, 1958.
\bibitem[M]{M} D.Mumford, {\em Rational equivalence for $0$-cycles on surfaces.}, J.Math Kyoto Univ. 9, 1968, 195-204.

\bibitem[R]{R} A.Roitman, {\em $\Gamma$-equivalence of zero dimensional cycles (Russian)}, Math. Sbornik. 86(128), 1971, 557-570.
\bibitem[R1]{R1}A.Roitman, {\em Rational equivalence of 0-cycles}, Math USSR Sbornik, 18, 1972, 571-588
\bibitem[R2]{R2} A.Roitman, {\em The torsion of the group of 0-cycles modulo rational equivalence}, Ann. of Math. (2), 111, 1980, no.3, 553-569
\bibitem[Ru1]{Ru1} K.Rubin, {\em Tate-Shafarevich group and L-functions of elliptic curves with complex multiplication}, Inventiones Mathematicae, 89(3), 527-559, 1987.
\bibitem[Sel]{Sel} E.Selmer, {\em On the diophantine equation $ax^3+by^3+cz^3=0$}, Acta Math. 85, 203-362, 1951.
\bibitem[Sil]{Sil} J.H.Silverman, {\em The Arithmetic of elliptic curves}, Springer, Berlin-Heidelberg-Newyork, 1986.
\bibitem[Sh1]{Sh1} I.R.Shafarevich, {\em The group of principal homogeneous algebraic manifolds}, Dokaldy Academii Nauk SSSR, (in Russian), 124, 42-43, 1959.
\bibitem[T]{T} J.Tate, {\em WC-groups over p-adic fields}, Seminaire Bourbaki, 1957-58, 13 Paris, Secretariat Mathematique.
\bibitem[SV]{SV} A.Suslin, V.Voevodsky, {\em Relative cycles and Chow sheaves}, Cycles, transfers, motivic homology theories, 10-86, Annals of Math studies.
\bibitem[Voi]{Voi}C.Voisin,{\em Symplectic invoultions of K$3$ surfaces act trivially on $\CH_0$}, Documenta Mathematicae 17, 851-860,2012.
\bibitem[Vo]{Vo} C.Voisin, {\em Complex algebraic geometry and Hodge theory II}, Cambridge studies of Mathematics, 2002.
\end{thebibliography}
\end{document}